\documentclass[11pt]{article}
\usepackage{graphicx}
\usepackage{amsfonts}
\usepackage{nath}
\usepackage{theorem}
\newtheorem{Lem}{Lemma}[section]

\newtheorem{The}[Lem]{Theorem}
\newtheorem{Prop}[Lem]{Proposition}

\newcommand{\qed}{\hbox{\rule{6pt}{6pt}}}

\begin{document}
\title{Bounds of logarithmic mean}
\author{Shigeru Furuichi$^1$\footnote{E-mail:furuichi@chs.nihon-u.ac.jp} \,\,and Kenjiro Yanagi$^2$\footnote{E-mail:yanagi@yamaguchi-u.ac.jp}\\
$^1${\small Department of Information Science,}\\
{\small College of Humanities and Sciences, Nihon University,}\\
{\small 3-25-40, Sakurajyousui, Setagaya-ku, Tokyo, 156-8550, Japan}\\
$^2${\small Division of Applied Mathematical Science,}\\
{\small Graduate School of Science and Engineering, Yamaguchi University,}\\
{\small 2-16-1, Tokiwadai, Ube City, 755-0811, Japan}}
\date{}
\maketitle
{\bf Abstract.} We give tight bounds for logarithmic mean. We also give new Frobenius norm inequalities for two positive semidefinite matrices. In addition, we give some matrix inequalities on matrix power mean.

\vspace{3mm}


\vspace{3mm}
\vspace{3mm}

\section{Introduction}
In this short paper, we study  the bounds of the logarithmic mean which is defined by
\begin{equation} \label{lm_def}
L(a,b) \equiv \frac{a-b}{\log a-\log b}= \int _0^1 a^{\nu} b^{1-\nu} d\nu, \quad a \neq b
\end{equation}
 for two positive numbers $a$ and $b$. (We conventionally define $L(a,b) = a$, if $a=b$.)
In the paper \cite{Lin},  the following relations were shown.
\begin{equation} \label{scalar_bounds_lm_01}
\sqrt {ab} \leq L(a,b)  \leq \left( \frac{a^{1/3}+b^{1/3}}{2}\right)^3, \quad a, b >0.
\end{equation}
We now have the following lemma.
\begin{Lem} \label{first_lemma}
For $a, b>0$, we have
\begin{equation} \label{scalar_upper_bounds_lm}
L(a,b) \leq \left( \frac{a^{1/3}+b^{1/3}}{2}\right)^3 \leq \frac{2}{3}\sqrt{ab}+ \frac{1}{3}\frac{a+b}{2}.
\end{equation}
\end{Lem}
{\it Proof:}
The second inequality of the inequalities (\ref{scalar_upper_bounds_lm}) can be proven easily. 
Indeed, we put,
$f(t) \equiv \frac{2}{3} t^3 + \frac{1}{6}(1+t^6)-\frac{1}{8}(1+t^2)^3 .$
Then we have $f'(t) <0 $ for $0<t<1$ and $f'(t) >0$ for $t>1$.
Thus we have $f(t) \geq f(1) =0$ for $t>0$.
\hfill \qed

The first inequality of the inequalities (\ref{scalar_upper_bounds_lm}) refines the inequality
\begin{equation}
 L(a,b)  \leq \frac{2}{3}\sqrt{ab}+ \frac{1}{3}\frac{a+b}{2},\quad a, b>0
\end{equation}
which is known as classical P\'olya inequality \cite{Zou, NP}.
 
Throughout this paper we use the notation $M(n,\mathbb{C})$
as the set of all $n \times n$ matrices on the complex field $\mathbb{C}$.
We also use  the notation $M_+(n,\mathbb{C})$ as  the set of all $n \times n$ positive semidefinite matrices. Here $A \in M_+(n,\mathbb{C})$ means we have $\langle \phi \vert A \vert \phi \rangle \geq 0$ for any vector
$\vert \phi \rangle \in \mathbb{C}^n$. For $A \in M(n,\mathbb{C})$,
the Frobenius norm (Hilbert-Schmidt norm) $\left\| \cdot \right\|_F$ is defined by
\begin{equation}
\left\| A \right\|_F \equiv \left( \sum_{i,j=1}^n \vert a_{ij}\vert ^2\right)^{1/2} =
\left(Tr\left[A^*A\right]\right)^{1/2}.
\end{equation}
In the paper \cite{Zou}, the following norm inequality was shown.
\begin{The} {\bf (\cite{Zou})} \label{Zou_norm_ineq_theorem}
For $A, B \in M_+(n,\mathbb{C})$, $X \in M(n,\mathbb{C})$ and Frobenius norm $\lVert \cdot \rVert_F$, we have
\begin{equation} \label{theorem2_Zou}
\lVert \int_0^1 {{A^{\nu}}X{B^{1 - \nu}}d\nu}  \rVert_F \le 
\frac{1}{3} \lVert {2{A^{1/2}}X{B^{1/2}} + \frac{{AX + XB}}{2}} \rVert_F.
\end{equation}
\end{The}
From Lemma \ref{first_lemma},
 we have the following proposition.
\begin{Prop}  \label{upper_bound_LM_prop}
For $A, B \in M_+(n,\mathbb{C})$, $X \in M(n,\mathbb{C})$ and Frobenius norm $\lVert \cdot \rVert_F$, we have
\begin{eqnarray*}
\lVert {\int_0^1 {{A^{\nu}}X{B^{1 - \nu}}d\nu} } \rVert_F &\le& 
\lVert {\frac{{AX + 3{A^{1/3}}X{B^{2/3}} + 3{A^{2/3}}X{B^{1/3}} + XB}}{8}} \rVert_F \\
&\le& \frac{1}{3}\lVert {2{A^{1/2}}X{B^{1/2}} + \frac{{AX + XB}}{2}}\rVert_F
\end{eqnarray*}
\end{Prop}
To the first author's best knowledge, the first inequality in Proposition \ref{upper_bound_LM_prop} was suggested in \cite{M_Lin}.

This proposition can be proven by the similar way to the proof of Theorem \ref{Zou_norm_ineq_theorem} (or the proof of Theorem \ref{Lower_Bound_main_theorem} which will be given in the next section) and this refines the inequality (\ref{theorem2_Zou}) shown in \cite{Zou}.


\section{Lower bound of logarithmic mean}
The following inequalities were given in \cite{HK}. Hiai and Kosaki gave the norm inequalities for Hilbert space operators in \cite{HK}. See also \cite{Hiai-Kosaki-book,Hiai2010}. Here we give them as a matrix setting to unify the description of this paper.
\begin{The} {\bf (\cite{HK})} \label{HK_Lower_Bound_theorem}
For $A, B \in M_+(n,\mathbb{C})$, $X \in M(n,\mathbb{C})$, $m \geq 1$ and
every  unitarily invariant norm $\ltriple|   \cdot  \rtriple|$, we have
$$
\ltriple|  {\int_0^1 {{A^{\nu}}X{B^{1 - \nu}}d\nu} }  \rtriple|\ge \frac{1}{m}
\ltriple|  {\sum\limits_{k = 1}^m {{A^{\frac{k}{{m + 1}}}}X{B^{\frac{{m + 1 - k}}{{m + 1}}}}} }  \rtriple| \ge \ltriple| {{A^{1/2}}X{B^{1/2}}}  \rtriple|.
$$
\end{The}
Frobenius norm is one of unitarily invariant norms. 
We give the refinement of the lower bound of the above first inequality
for Frobenius norm. That is, we have the following inequalities.

\begin{The} \label{Lower_Bound_main_theorem}
For $A, B \in M_+(n,\mathbb{C})$, $X \in M(n,\mathbb{C})$, $m \geq 1$ and Frobenius norm $\lVert \cdot \rVert_F$, we have
\begin{eqnarray}
\lVert {\int_0^1 {{A^{\nu}}X{B^{1 - \nu}}d\nu} } \rVert_F &\ge& 
\frac{1}{m}\lVert {\sum\limits_{k = 1}^m {{A^{\frac{{2k - 1}}{{2m}}}}X{B^{\frac{{2m - \left( {2k - 1} \right)}}{{2m}}}}} } \rVert_F \nonumber \\
&\ge& 
\frac{1}{m}\lVert {\sum\limits_{k = 1}^m {{A^{\frac{k}{{m + 1}}}}X{B^{\frac{{m + 1 - k}}{{m + 1}}}}} }\rVert_F 
\ge \lVert {{A^{1/2}}X{B^{1/2}}} \rVert_F. \label{Lower_Bound_main_theorem_ineq00}
\end{eqnarray}
\end{The}

To prove Theorem \ref{Lower_Bound_main_theorem}, we need a few lemmas.
\begin{Lem}\label{Lower_Bound_lemma01}
Let  $u,v,w$  be nonnegative integers such that $w \geq u $ and let $x$ be a positive real number. Then we have
\begin{equation}
x^u (1-x^v) + x^w(x^v -1) \geq 0.
\end{equation}
\end{Lem}
{\it Proof :}
It is trivial for the case $x=1$ or $v=0$. We prove for the case $x \neq 1$ and $v \neq 0$.
In addition, for the case that $u=w$, the equality holds. Thus we may assume
$w > u$ and $v \geq 1$. Then the lemma can be proven by the following way.
\begin{eqnarray*}
&&x^u ( 1 - x^v ) + x^w ( x^v - 1 ) = ( x^v - 1 )( x^w - x^u) 
= x^u ( x^v - 1 )( x^{w - u} - 1)\\
&& = x^u ( x - 1)^2 ( x^{v - 1} + x^{v - 2} +  \cdots  + 1)
( x^{w - u - 1} + x^{w - u - 2} +  \cdots  + 1) \ge 0.
\end{eqnarray*}
\hfill \qed

\begin{Lem}\label{Lower_Bound_lemma02}
For a positive real number $x$ and a natural number $m$, we have
\begin{equation}\label{Lower_Bound_lemma02_ineq01}
\sum\limits_{k = 1}^m {{x^{\left( {2k - 1} \right)\left( {m + 1} \right)}}}  \ge \sum\limits_{k = 1}^m {{x^{2km}}} 
\end{equation}
\end{Lem}
{\it Proof :}
For the case $x=1$, the equality holds.  So we prove this lemma for $x \neq 1$.
If $m$ is an odd number, then we have
$m = 2 \lfloor m /2 \rfloor + 1$. Since we then have
$m+1 = 2 (\lfloor m/2 \rfloor +1)$ and $2 (\lfloor m/2 \rfloor +1)-1 =m$,
we have
$$
\{2 (\lfloor m/2 \rfloor +1)-1\}(m+1) = 2 (\lfloor m/2 \rfloor +1) m.
$$
If we put $\tilde{k} = \lfloor m/2 \rfloor +1$, then the above means
$(2\tilde{k}-1)(m+1)=2\tilde{k}m$. Then the difference of the $\tilde{k}$-th term of the both sides in the inequality (\ref{Lower_Bound_lemma02_ineq01})
is equal to $0$.
For the case that $m$ is an even number, it never happens that the difference of the $\tilde{k}$-th term of the both sides in the inequality (\ref{Lower_Bound_lemma02_ineq01}) is equal to $0$.
Therefore we have
\begin{equation}\label{Lower_Bound_lemma02_ineq02}
\sum\limits_{k = 1}^m x^{( 2k - 1)( m + 1)}  - \sum\limits_{k = 1}^m x^{2km} =
 \sum_{l=1}^{\lfloor m/2 \rfloor} \{x^{a_l}(1-x^{b_l}) +x^{c_l}(x^{b_l} -1)\},
 \end{equation}
where 
$a_l=(2l-1)(m+1)$, $b_l=m-(2l-1)$ and $c_l=2\{m-(l-1)\}m$, for $l=1,2,\cdots,\lfloor m/2 \rfloor$.
Here we have $c_l-a_l = \{ m-(2l-1) \}(2m+1) \ge 0$ whenever $b_l \geq 0$.
By Lemma \ref{Lower_Bound_lemma01}, if  $b_l \geq 0$, then we have
$x^{a_l}(1-x^{b_l}) +x^{c_l}(x^{b_l} -1) \geq 0$. Thus the proof of this lemma was completed.

\hfill \qed

If we put $t=x^{2m(m+1)} >0$, then we have 
\[\sum\limits_{k = 1}^m {{t^{\frac{{2k - 1}}{{2m}}}}}  \ge \sum\limits_{k = 1}^m {{t^{\frac{k}{{m + 1}}}}}, \,\,\,\,\left( {t > 0,m \in \mathbb{N}} \right),\]
which implies
\begin{equation} \label{Lower_Bound_lemma02_ineq03}
\sum\limits_{k = 1}^m {{a^{\frac{{2k - 1}}{{2m}}}}{b^{\frac{{2m - \left( {2k - 1} \right)}}{{2m}}}}}  \ge \sum\limits_{k = 1}^m {{a^{\frac{k}{{m + 1}}}}{b^{\frac{{m + 1 - k}}{{m + 1}}}}}, \,\,\,\,\left( {a,b > 0,m \in \mathbb{N}} \right).
\end{equation}

We then have the following lemma.
\begin{Lem} \label{Lower_Bound_lemma03}
For $a,b >0$ and $m \geq 1$, we have
\begin{equation}\label{Lower_Bound_lemma03_ineq01}
\hspace*{-10mm}L(a,b) \ge \frac{1}{m}\sum\limits_{k = 1}^m {{a^{\frac{{2k - 1}}{{2m}}}}{b^{\frac{{2m - \left( {2k - 1} \right)}}{{2m}}}}}  \ge \frac{1}{m}\sum\limits_{k = 1}^m {{a^{\frac{k}{{m + 1}}}}{b^{\frac{{m + 1 - k}}{{m + 1}}}}}\ge \sqrt {ab}. 
\end{equation}
\end{Lem}
{\it Proof :}
The second inequality follows by the inequality (\ref{Lower_Bound_lemma02_ineq03}).
We use the famous inequality $\frac{x-1}{\log x} \geq \sqrt{x} $ for $x>0$.
We put $x = t^{1/m}$ in this inequality. Then we have for $t >0$
\begin{equation}  \label{appendix_sum01}
\frac{{t - 1}}{{\log t}} \ge \frac{{{t^{\frac{{2m - 1}}{{2m}}}} + {t^{\frac{{2m - 3}}{{2m}}}} +  \cdots  + {t^{\frac{1}{{2m}}}}}}{m},
\end{equation}
which implies the first inequality.
The third inequality can be proven by the use of the arithmetic mean - geometric mean inequality. Thus the proof of this lemma was completed.

\hfill \qed

We give some basic properties of the right hand side of the  inequality (\ref{appendix_sum01}) in Appendix.
\\

{\it Proof of Theorem \ref{Lower_Bound_main_theorem}: }
It is known that the Frobenius norm inequality immediately follows from the corresponding scalar inequality \cite{Hiai-Kosaki-book}. However, we give here an elementary proof for the convenience of the readers. 
Let ${U^*}AU = {D_1} = diag\left\{ {{\alpha _1}, \cdots ,{\alpha _n}} \right\}$
 and ${V^*}BV = {D_2} = diag\left\{ {{\beta _1}, \cdots ,{\beta _n}} \right\}$.
Then ${\alpha _1}, \cdots ,{\alpha _n} \ge 0$ and ${\beta _1}, \cdots ,{\beta _n} \ge 0$. We put ${U^*}XV \equiv Y = \left( {{y_{ij}}} \right)$. By the first inequality of the inequalities (\ref{Lower_Bound_lemma03_ineq01}), we have
\begin{eqnarray*}
&& \hspace*{-15mm} \frac{1}{m} \lVert \sum\limits_{k = 1}^m 
A^{\frac{2k - 1}{2m}} X  B^{\frac{ 2m - ( 2k - 1)}{2m}}  \rVert_F \\
&&\hspace*{-15mm}= \frac{1}{m}\lVert \sum\limits_{k = 1}^m 
U D_1^{\frac{2k - 1}{2m}} U^* XV D_2^{\frac{2m -( 2k - 1)}{2m}}V^* \rVert_F \\
&&\hspace*{-15mm}= \frac{1}{m} \lVert \sum\limits_{k = 1}^m 
D_1^{\frac{2k - 1}{2m}}Y D_2^{\frac{2m -( 2k - 1)}{2m}}  \rVert_F \\
&&\hspace*{-15mm}=\frac{1}{m} \lVert ( \begin{array}{*{3}{c}}
\sum\limits_{k = 1}^m \alpha _1^{\frac{2k - 1}{2m}}
y_{11}\beta _1^{\frac{2m - ( 2k - 1)}{2m}} & \ldots &\sum\limits_{k = 1}^m \alpha _1^{\frac{2k - 1}{2m}}y_{1n}\beta _n^{\frac{2m - ( {2k - 1} )}{2m}} \\
 \vdots & \ddots & \vdots \\
\sum\limits_{k = 1}^m \alpha _n^{\frac{2k - 1}{2m}}y_{n1}\beta _1^{\frac{2m - ( 2k - 1 )}{2m}} & \cdots &\sum\limits_{k = 1}^m \alpha _n^{\frac{2k - 1}{2m}}y_{nn}\beta _n^{\frac{2m - ( 2k - 1 )}{2m}} \end{array} ) \rVert_F \\
&&\hspace*{-15mm}= \{ \sum\limits_{i,j = 1}^n 
( \frac{1}{m}\sum\limits_{k = 1}^m \alpha _i^{\frac{2k - 1}{2m}}
\beta _j^{\frac{2m - ( 2k - 1 )}{2m}}  )^2 \lvert y_{ij} \rvert^2  \}^{1/2} \\
&& \hspace*{-15mm} \leq \{ \sum\limits_{i,j = 1}^n ( \int_0^1 
\alpha _i^{\nu}\beta _j^{1 - \nu}d\nu  )^2 \lvert y_{ij} \rvert^2 \}^{1/2}
 = \lVert \int_0^1 D_1^{\nu}YD_2^{1 - \nu}d\nu  \rVert_F
 = \lVert U \int_0^1 D_1^{\nu}YD_2^{1 - \nu}d\nu  V^*\rVert_F \\
&&\hspace*{-15mm}=\lVert  \int_0^1 UD_1^{\nu} U^*XVD_2^{1 - \nu}V^* d\nu  \rVert_F=\lVert  \int_0^1A^{\nu} XB^{1 - \nu} d\nu  \rVert_F.
\end{eqnarray*}
Applying the inequality (\ref{Lower_Bound_lemma02_ineq03}), we have
the second inequality of (\ref{Lower_Bound_main_theorem_ineq00}) by the similar way. The third inequality holds due to Theorem \ref{HK_Lower_Bound_theorem} (or the third inequality of the inequalities (\ref{Lower_Bound_lemma03_ineq01})).

\hfill \qed


\section{Upper bound of logarithmic mean}
In the paper \cite{HK}, the following norm inequalities were also given for Hilbert space operators. Here we give them for matrices, as we mentioned in the beginning of Section 2.

\begin{The} {\bf (\cite{HK})} \label{HK_Upper_Bound_theorem}
For $A, B \in M_+(n,\mathbb{C})$, $X \in M(n,\mathbb{C})$, $m \geq 2$ and
every  unitarily invariant norm $\ltriple|   \cdot  \rtriple|$, we have
$$
\ltriple|  {\int_0^1 {{A^{\nu}}X{B^{1 - \nu}}d\nu} }  \rtriple|
\le \frac{1}{m}
\ltriple|  {\sum\limits_{k = 0}^{m-1} {{A^{\frac{k}{{m - 1}}}}X{B^{\frac{{m - 1 - k}}{{m - 1}}}}} }  \rtriple| 
\le \frac{1}{2} \ltriple| AX +XB  \rtriple|.
$$
\end{The}

We also give an improved upper bound of the logarithmic mean on the above Theorem \ref{HK_Upper_Bound_theorem}, only for the Frobenius norm. Namely, we can prove the following inequalities by the similar way to the proof of Theorem \ref{Lower_Bound_main_theorem}, by the use of scalar inequalities 
which will be given in Lemma \ref{Upper_Bound_lemma01}.

\begin{The} \label{Upper_Bound_main_theorem}
For $A, B \in M_+(n,\mathbb{C})$, $X \in M(n,\mathbb{C})$, $m \geq 2$ and Frobenius norm $\lVert \cdot \rVert_F$, we have
\begin{eqnarray}
\lVert \int_0^1 A^{\nu} X B^{1 - \nu} d\nu  \rVert_F &\le& 
\frac{1}{m} \lVert \sum\limits_{k = 0}^m A^{\frac{k}{m}} X B^{\frac{m - k}{m}}  - \frac{1}{2} (  AX + XB ) \rVert_F \nonumber \\
&\le&  \frac{1}{m}
\lVert  {\sum\limits_{k = 0}^{m-1} {{A^{\frac{k}{{m - 1}}}}X{B^{\frac{{m - 1 - k}}{{m - 1}}}}} } \rVert_F \leq \frac{1}{2} \lVert AX + BX \rVert_F. 
 \label{Upper_Bound_main_theorem_ineq00}
\end{eqnarray}
\end{The}

To prove Theorem \ref{Upper_Bound_main_theorem} we need to prove the following lemmas.

\begin{Lem}  \label{conjectured_lemma}
For $x>0$ and $m \geq 2$, we have
\begin{equation} \label{con_lemma_eq00}
\sum\limits_{k = 1}^{m - 1} ( x^{km} - x^{k ( m - 1)} )  \ge 
\frac{ x^{m ( m - 1 )} - 1}{2},\quad ( x > 0 ).
\end{equation}
\end{Lem}
{\it Proof:}
For $m \geq 2$, we calculate
\begin{eqnarray}
&& \sum\limits_{k = 1}^{m - 1} ( x^{km} - x^{k ( m - 1)} )  - \frac{ x^{m ( m - 1 )} - 1}{2} \nonumber \\
&& =\frac{1}{2} (x^m +1 )\{ (x^m)^{m-2} +\cdots +x^m+1\} -x^{m-1} \{ (x^{m-1})^{m-2} + \cdots +x^{m-1} +1\} \nonumber \\
&& \geq x^{m/2}\{ (x^m)^{m-2} +\cdots +x^m+1\} -x^{m-1} \{ (x^{m-1})^{m-2} + \cdots +x^{m-1} +1\}.  \label{con_lemma_eq01}
\end{eqnarray}
Here we put $y \equiv x^{1/2} >0$, then we have
\begin{eqnarray}
(\ref{con_lemma_eq01}) &=& y^m\{ (y^{2m})^{m-2} + \cdots +y^{2m} +1\} -y^{2(m-1)} \{ (y^{2(m-1)})^{m-2} +\cdots +y^{2(m-1)} +1\} \nonumber \\
&=& y^m (1-y^{m-2}) + y^{3m}(1-y^{m-4}) +y^{5m} (1-y^{m-6})+ \cdots \nonumber \\
&&+y^{2(m-1)(m-3)}(y^{m-6} -1)+y^{2(m-1)(m-2)} (y^{m-4}-1)+y^{2(m-1)^2}(y^{m-2} -1) \nonumber \\
&=& \sum_{l=1}^{\left\lfloor {m/2} \right\rfloor} \{ y^{p_l} (1-y^{r_l}) + y^{q_l} (y^{r_l
}-1)\} \geq 0,
\end{eqnarray}
where $p_l = (2l-1)m$, $q_l=2(m-1)(m-l)$ and $r_l = m-2l$.
The last inequality follows from Lemma \ref{Lower_Bound_lemma01},
because we have $q_l-p_l = (2m-1)(m-2l) \geq $ whenever $r_l \geq 0$.

\hfill \qed

\begin{Lem} \label{Upper_Bound_lemma01}
For $a,b >0$ and $m \geq 2$, we have
\begin{equation} \label{Upper_Bound_lemma01_ineq01}
\hspace*{-15mm} L(a,b)  \le 
\frac{1}{m} \{ \sum\limits_{k = 0}^m a^{k/m} b^{( m - k )/m}  - 
( \frac{a + b}{2} ) \} \le \frac{1}{m} \sum_{k=0}^{m-1}a^{k/(m-1)}b^{(m-1-k)/(m-1)} \leq \frac{a + b}{2}.
\end{equation}
\end{Lem}
{\it Proof :}
To prove the first inequality, we have only to prove the following inequality
\begin{equation} \label{Upper_Bound_lemma01_ineq02}
\frac{t - 1}{\log t} \le \frac{1}{m}\{ ( t + t^{( m - 1)/m} + t^{( m - 2)/m} +  \cdots  + t^{1/m} + 1 ) - ( \frac{t + 1}{2} ) \},\quad (t > 0).
\end{equation}
The inequality  (\ref{Upper_Bound_lemma01_ineq02})
can be proven by putting $x \equiv t^{1/m} >0 $ in a famous inequality $\frac{x-1}{\log x} \leq \frac{x+1}{2}$ for $x >0$. 

To prove the second inequality of the inequalities (\ref{Upper_Bound_lemma01_ineq01}), it is sufficient to prove 
the inequality (\ref{con_lemma_eq00}) which holds from Lemma \ref{conjectured_lemma}. We obtain actually the second inequality of the inequalities (\ref{Upper_Bound_lemma01_ineq01}) by putting $t= x^{m(m-1)} >0 $
in the inequality (\ref{con_lemma_eq00}), and then putting $t = a/b$.

To prove the third inequality of the inequalities (\ref{Upper_Bound_lemma01_ineq01}),  it is sufficient to prove 
the following inequality
\begin{equation} \label{con_lemma_induction01}
t^{m-1} + \cdots +t +1 \leq \frac{m(t^{m-1}+1)}{2}, \quad (t>0).
\end{equation}
This inequality can be proven by the induction on $m$.
Indeed, we assume that the inequality (\ref{con_lemma_induction01}) holds for some $m$. Then we add $t^m >0$ to both sides to he inequality (\ref{con_lemma_induction01}). Then we have
$$
t^m +t^{m-1} + \cdots t+1 \leq \frac{m(t^{m-1}+1)}{2} +t^m.
$$
Therefore we have only to prove the inequality
$$
\frac{m(t^{m-1}+1)}{2} +t^m \leq \frac{(m+1)(t^m+1)}{2},\quad (t>0)
$$
which is equivalent to the inequality
$$
(m-1)t^m -m t^{m-1} +1 \geq 0, \quad (t>0). 
$$
We put $f_m(t)\equiv (m-1)t^m -m t^{m-1} +1$. Then we can prove $f_m(t) \geq f_m(1) =0$ by elementary calculations.
Thus the inequality (\ref{con_lemma_induction01}) holds for $m+1$.

\hfill \qed

We give some basic properties of the right hand side of the  inequality (\ref{Upper_Bound_lemma01_ineq02}) in Appendix.


\section{Matrix inequalities on geometric mean}
Using Lemma \ref{first_lemma}, Lemma \ref{Lower_Bound_lemma03} and Lemma \ref{Upper_Bound_lemma01}, we have following Proposition \ref{sec4_theorem01},  Proposition \ref{sec4_theorem02} and Proposition \ref{sec4_theorem03}, respectively.

\begin{Prop} \label{sec4_theorem01}
For $A, B \in M_+(n,\mathbb{C})$, we have
$$
\int_0^1 A\#_{\nu}Bd \nu \leq \frac{1}{4} \{ \frac{A+B}{2} + \frac{3}{2}(A\#_{2/3}B+A\#_{1/3}B)\} \leq \frac{1}{2} (\frac{A+B}{2}+2A\#_{1/2}B),
$$
where $A\#_{\nu}B \equiv A^{1/2}(A^{-1/2}BA^{-1/2})^{\nu}A^{1/2}$, $(0 \leq \nu \leq 1)$ is $\nu$-weighted geometric mean introduced in \cite{KA}.
\end{Prop}

\begin{Prop} \label{sec4_theorem02}
For $A, B \in M_+(n,\mathbb{C})$  and $m \geq 1$, we have
$$
\int_0^1 A\#_{\nu}Bd \nu \geq \frac{1}{m}\sum_{k=1}^m A\#_{(2k-1)/2m}B \geq \frac{1}{m} \sum_{k=1}^mA\#_{k/(m+1)}B \geq A\#_{1/2}B.
$$
\end{Prop}

\begin{Prop} \label{sec4_theorem03}
For $A, B \in M_+(n,\mathbb{C})$  and $m \geq 2$, we have
$$
\int_0^1 A\#_{\nu}Bd \nu \leq \frac{1}{m} (\sum_{k=0}^mA\#_{k/m}B -\frac{A+B}{2})\leq
\frac{1}{m}\sum_{k=0}^{m-1}A\#_{k/(m-1)}B \leq
 \frac{A+B}{2}.
$$
\end{Prop}

We give the proof of Proposition \ref{sec4_theorem03}.
Proposition \ref{sec4_theorem01} and Proposition \ref{sec4_theorem02} are also proven by the similar way, using  Lemma \ref{first_lemma} and Lemma \ref{Lower_Bound_lemma03}. In addition, by using the notion of the representing function $f_m(x)=1 m x$ for operator mean $m$, it is well-known \cite{KA} that $f_m(x) \leq f_n(x)$ holds for $x>0$ if and only if $A m B \leq A n B$ holds for all positive operators $A$ and $B$. However we give an elementary proof for the convenience of the readers.

{\it Proof of Proposition \ref{sec4_theorem03}:}
Since $T \equiv A^{-1/2}BA^{-1/2} \geq 0$, there exists a unitary matrix $U$
such that $U^*TU = D \equiv diag \{\lambda_1,\cdots,\lambda_n\}$. Then
$\lambda_1,\cdots,\lambda_n \geq 0$. From Lemma \ref{Upper_Bound_lemma01}, for $i=1,\cdots,n$, we have
$$
\int_0^1 \lambda_i^{\nu}d\nu \leq \frac{1}{m}\{\sum_{k=0}^m \lambda_i^{k/m} -(\frac{\lambda_i +1}{2})\} \leq \frac{1}{m}\sum_{k=0}^{m-1} \lambda_i^{k/(m-1)}  \leq \frac{\lambda_i +1}{2}.
$$
Thus we have
$$
\int_0^1 D^{\nu}d\nu \leq \frac{1}{m}\{\sum_{k=0}^m D^{k/m} -(\frac{D +I}{2})\} \leq \frac{1}{m}\sum_{k=0}^{m-1} D^{k/(m-1)} \leq \frac{D +I}{2}.
$$
Multiplying $U$ and $U^*$ to both sides, we have
$$
\int_0^1 T^{\nu}d\nu \leq \frac{1}{m}\{\sum_{k=0}^m T^{k/m} -(\frac{T +I}{2})\} \leq \frac{1}{m}\sum_{k=0}^{m-1} T^{k/(m-1)} \leq \frac{T +I}{2}.
$$
Inserting $T \equiv A^{-1/2}BA^{-1/2}$ and then multiplying two $A^{1/2}$ to all sides from both sides, we obtain the result.

\hfill \qed

Closing this section, we give another matrix inequalities by the use of the another lower bound of the logarithmic mean.
As an another lower bound of the logarithmic mean, the following inequalities are known.
\begin{equation}\label{conclusion_ineq01}
\frac{t-1}{\log t} \geq \frac{t+t^{1/3}}{1+t^{1/3}} \geq \sqrt {t}, \quad t>0,\,\, t\neq 1.
\end{equation}
The  proofs of the above inequalities are not so difficult, (they can be done by putting $x = t^{1/3} >0$ and $x = t^{1/6} >0$)
here we omit to write them. From the inequalities (\ref{conclusion_ineq01}),  we have
\begin{equation}\label{conclusion_ineq02}
\int_0^1 t^{\nu} d\nu \geq t^{2/3} -t^{1/3}+2(t^{-1/3}+1)^{-1} \geq \sqrt{t},\quad t>0,t\neq 1.
\end{equation}
The inequalities (\ref{conclusion_ineq02}) imply the following result by the similar way to the proof of Proposition \ref{sec4_theorem03}.
\begin{Prop} \label{sec4_theorem04}
For $A, B \in M_+(n,\mathbb{C})$, we have
$$
\int_0^1 A\#_{\nu}Bd \nu \geq A\#_{2/3}B-A\#_{1/3}B+2
\{ (A\#_{1/3}B)^{-1}+A^{-1}\}^{-1} \geq A\#_{1/2}B.
$$
\end{Prop}


\section{Comments}
Proposition \ref{appendix_prop2} given in Appendix shows that our upper bound is tighter than the standard upper bound for the case $t>1$ and our lower bound is tighter than the standard lower bound for any $t > 0$. 
In addition, our lower bound $\alpha_m(t)$ of the logarithmic mean $L(t,1)$ is tighter than the lower bound $\sqrt{t}$ given by T.-P.Lin in \cite{Lin}, for $m \geq 1$. However, it may be difficult problem to find the minimum $m \in \mathbb{N}$  such that 
$\beta_m(t) \leq (\frac{t^{1/3}+1}{2})^3$ for any $t >0$.
The right hand side of the above inequality is the upper bound given  by T.-P.Lin in \cite{Lin}. (See the inequalities (\ref{scalar_bounds_lm_01}).)

\section*{Competing interests}
The authors declare that they have no competing interests.

\section*{Authors’ contributions}
The work presented here was carried out in collaboration between all authors. The study was initiated by S.F. and the manuscript was written by S.F.
The author S.F. also played the role of the corresponding author. 
The proof of Lemma \ref{first_lemma} and the first equality of Eq.(\ref{con_lemma_eq01}) in the proof of Lemma \ref{conjectured_lemma} were given by the author K.Y. With the exception of them, the proofs of all results were given by the author S.F.
All authors have contributed to, checked and approved the manuscript.

\section*{Ackowledgements}
The author (S. F.) was partially supported by JSPS KAKENHI Grant Number 24540146. The author (K.Y.) was also partially supported by JSPS KAKENHI Grant Number 23540208.

\section*{Appendix}
Here we note some basic properties of the following scalar sums
$$
\alpha_m (t)\equiv \frac{1}{m} \sum_{k=1}^m t^{(2k-1)/(2m)}, \quad  \quad 
\beta_m (t)\equiv \frac{1}{m} \left(\sum_{k=0}^m t^{k/m} -\frac{t+1}{2}\right)
$$
for $t > 0$.
\begin{Prop}
For any $t > 0$, we have following properties:
\begin{itemize}
\item[(i)] $\alpha_m (t) \leq \alpha_{m+1}(t).$
\item[(ii)] $\beta_{m+1}(t) \leq \beta_m(t).$
\item[(iii)] $\alpha_m(t)$ and $\beta_m(t)$ converges to $L(t,1)$ as $m \rightarrow \infty $. In addition, we have $\alpha_m(t) \leq \beta_m(t)$.
\end{itemize}
\end{Prop}

{\it Proof:}
We prove (i)-(iii) for $t \neq 1$, since it is trivial for the case $t=1$.
\begin{itemize}
\item[(i)]
Since 
$$
\alpha_m(t) = \frac{ t^{1/(2m)} (t-1) }{m(t^{1/m}-1)},
$$
for $t\neq 1$, we prove 
$$
\frac{\alpha_{m+1}(t)}{\alpha_m(t)} = \frac{mt^{1/(2(m+1))} (t^{1/m}-1)}{(m+1)t^{1/(2m)} (t^{1/(m+1)}-1)} > 1,
$$
for $t >0$ and $t \neq 1$.
We first prove the case $t > 1$. Then we put $s \equiv t^{1/(2m(m+1))}$ and
$$
f_m(s) \equiv m (s^{2m+2}-1)-(m+1)(s^{2m+1}-s),\quad (s>1).
$$
By elementary calculations, we have $f_m(s) > f_m(1) =0$, which implies
$$
ms^m(s^{2m+2} -1) > (m+1)s^{m+1}(s^{2m} -1).
$$
We can prove similarly
$$
ms^m(s^{2m+2} -1) < (m+1)s^{m+1}(s^{2m} -1),
$$
for the case $0<s<1$.
\item[(ii)] Since 
$$
\beta_m(t) = \frac{(t^{1/m}+1)(t-1) }{2m(t^{1/m}-1)},
$$
for $t\neq 1$, we prove
$$
\frac{\beta_{m+1}(t)}{\beta_m(t)} = \frac{m(t^{1/(m+1)}+1)(t^{1/m}-1)}{(m+1)(t^{1/m}+1)(t^{1/(m+1)}-1)} < 1
$$
for $t >0$ and $t \neq 1$.
We first prove the case $t >1$. Then we put $r \equiv t^{1/(m(m+1))}$ and
$$
g_m(r) \equiv r^{2m+1} -(2m+1)r^{m+1} +(2m+1)r^{m} -1,\quad (r >1).
$$
Since
$$
g'_m(r) = (2m+1)r^{m-1} (r^{m+1}-(m+1)r+m) >0, 
$$
we have $g_m(r) > g_m(1) =0$, which implies
$$
m(r^{m}+1)(r^{m+1}-1) < (m+1)(r^{m+1}+1)(r^m-1).
$$
We can prove similarly
$$
m(r^{m}+1)(r^{m+1}-1) > (m+1)(r^{m+1}+1)(r^m-1)
$$
for the case $0< r < 1$.
\item[(iii)] Since we have
$$\lim_{m \rightarrow \infty} m(t^{1/m}-1) =\lim_{p \rightarrow 0} \frac{t^p -1}{p} =\log t$$
 for $t >0$ and $t \neq 1$, we have 
$$
\lim_{m \rightarrow \infty} \alpha_m(t) = L(t,1),\quad \lim_{m \rightarrow \infty} \beta_m(t) = L(t,1).
$$
\end{itemize}
The arithmetic-geometric mean inequality proves $\alpha_m(t) < \beta_m(t)$, for $t >0$ and $t\neq 1$.

\hfill \qed
\\

As standard bounds of Riemann sum for the integral $\int_0^1 t^{\nu} d\nu$,
we have
$$
\gamma_m(t) < \int_0^1 t^{\nu} d\nu < \delta_m(t), \quad (0<t<1)
$$
and
$$
\delta_m(t) < \int_0^1 t^{\nu} d\nu < \gamma_m(t), \quad (t>1)
$$
where
$$
\gamma_m(t) \equiv \frac{1}{m} \sum_{k=1}^{m} t^{k/m}, \quad  \delta_m(t) \equiv \frac{1}{m} \sum_{k=0}^{m-1} t^{k/m}. 
$$
Then we have the following relations.
\begin{Prop} \label{appendix_prop2}
\begin{itemize}
\item[(i)] For $0<t<1$, we have $\alpha_m(t) > \gamma_m(t)$ and $\beta_m(t) < \delta_m(t)$.
\item[(ii)] For $t>1$, we have $\alpha_m(t) > \delta_m(t)$ and $\beta_m(t) < \gamma_m(t)$.
\end{itemize}
\end{Prop}

{\it Proof:}
\begin{itemize}
\item[(i)] For the case $0<t<1$, the following calculations shows assertion: 
$$\alpha_m(t)-  \gamma_m(t)=\frac{(t-1)(t^{1/(2m)}-t^{1/m})}{m(t^{1/m}-1)} >0, \quad 
\beta_m(t) - \delta_m(t) =\frac{t-1}{2m} <0.
$$
\item[(ii)]  For the case $t>1$, the following calculations shows assertion: 
$$\alpha_m(t)-  \delta_m(t)=\frac{t+t^{1/(2m)}}{m(t^{1/(2m)}+1)} >0, \quad 
\beta_m(t) - \gamma_m(t) =\frac{1-t}{2m} <0.
$$
\end{itemize}

\hfill \qed

\end{document}